\documentclass{article}

\def\l{\lambda}

\newcommand{\ds}{\displaystyle}

\newtheorem{theorem}{Theorem}
\newtheorem{example}{Example}

\newcommand{\be}{\begin{equation}}
\newcommand{\ee}{\end{equation}}

\newcommand{\D}{\displaystyle}

\newcommand{\beq}{\begin{eqnarray}}
\newcommand{\eeq}{\end{eqnarray}}
\newcommand{\nbeq}{\begin{eqnarray*}}
\newcommand{\neeq}{\end{eqnarray*}}
\def\st{\stackrel}
\def\CB{{\cal B}}
\def\CM{{\cal M}}

\def\CR{{\cal R}}
\def\CV{{\cal V}}
\def\i{\infty}

\begin{document}

 \title{\bf Revisiting Offspring Maxima in Branching Processes}
\author{George P. Yanev \\
Department of Mathematics and Statistics \\
University of South Florida \\
Tampa, Florida 33620 \\
e-mail: gyanev@cas.usf.edu}
\date{\empty}

\maketitle

\begin{abstract}
We present a progress report for studies on maxima related to
offspring in branching processes. We summarize and discuss the
findings on the subject that appeared in the last ten years. Some
of the results are refined and illustrated with new examples.
\end{abstract}

\section{Introduction}
There is a significant amount of research in the theory of
branching processes devoted to extreme value problems concerning
different population characteristics. The history of such studies
goes back to the works in 50-ies by Zolotarev \cite{Zol54} and
Urbanik \cite{Urb56} (see also \cite{Har63}) who considered the
maximum generation size. Our goal here is to summarize and discuss
results on maxima related to the offspring. Papers directly
addressing this area of study have begun to appear in the last ten
years (though see "hero mothers" example in \cite{JagNer84}.)

Let $\CM_n$ denote the maximum offspring size of all individuals
living in the $(n-1)$-st generation of a branching process. This
is a maximum of random number of independent and identically
distributed (i.i.d.) integer-valued random variables, where the
random index is the population size of the process. $\CM_n$ has
two characteristic features: (i) the i.i.d. random variables are
integer-valued and (ii) the distribution of the random index is
connected to the distribution of the terms involved through the
branching mechanism. These two characteristics distinguish the
subject matter maxima among those studied in the general extreme
value theory.

The study of the sequence \ $\{\CM_n\}$\ might be motivated in
different ways. It provides a fertility measure characterizing the
most prolific individual in one generation. It also measures the
maximum litter (or family) size.  In the branching tree context,
it is the maximum degree of a vertex. The asymptotic behavior of
$\CM_n$\ gives us some information about the influence of the
largest families on the size and survival of the entire
population.

The paper is organized as follows. Next section deals with maxima
in simple branching processes with or without immigration. In
Section~3 we derive results about maxima of a triangular array of
zero-inflated geometric variables. Later we apply these to
branching processes with varying geometric environments. Section 4
begins with limit theorems for the max-domain of attraction of
bivariate geometric variables. Then we discuss one application to
branching processes with promiscuous matting. The final section
considers a different construction in which a random score (a
continuous random variable) is associated with each individual in
a simple branching process. We present briefly limiting results
for the score's order statistics. In the end of the section, we
give an extension to two-type processes.

\section{Maximum family size in simple branching processes}

Define a Bienaym\'{e}--Galton--Watson (BGW) branching process and
its $n$-th generation maximum family size by $Z_0=1$;
\[ Z_n=\sum_{i=1}^{Z_{n-1}}X_i(n)\quad \mbox{and}\quad \CM_n=
\max_{ 1\le i\le Z_{n-1}}X_i(n) \quad (n=1,2,\ldots),
\]
respectively, where the offspring variables $X_i(n)$ are i.i.d.
nonnegative and integer-valued.

Along with the BGW process $\{Z_n\}$, we consider the process with
immigration $\{Z^{im}_n\}$ and its offspring maximum
\[ Z^{im}_n=\sum_{i=1}^{Z^{im}_{n-1}}X_i(n)+ Y_n \quad \mbox{and}
\quad \CM^{im}_n=\max_{1\le i\le Z^{im}_{n-1}}X_i(n) \quad
(n=1,2,\ldots),
\]
respectively, where $\{Y_n, \  n=1,2,...\}$ are independent of the
offspring variables, i.i.d. and integer-valued non-negative random
variables.

Finally, let us modify the immigration component such that
immigrants may enter the $n$-th generation only
 if the $(n-1)$-st generation size is zero. Thus, we have the Foster-Pakes
 process and its offspring maximum
\[ Z^{0}_n=\sum_{i=1}^{Z^0_{n-1}}X_i(n)+ I_{\D \{Z^{0}_{n-1}=0\}}Y_n
\quad \mbox{and} \quad \CM^{0}_n=\max_{1\le i \le
Z^0_{n-1}}X_i(n)\quad (n=1,2,\ldots),
\]
where $I_{A} $ stands for the indicator of $A$.

 Denote by $F(x)=P(X_i(n)\leq x)$\ the common distribution function of
the offspring variables with mean \ $0<m<\infty$\ and variance \
$0<\sigma^2 \leq \infty$. In this section, we deal with the
subcritical $(m<1)$, critical $(m=1)$, and supercritical $(m>1)$
processes separately.

\subsection{Subcritical processes}
Let $\hat{\CM}_n$ denote the maximum family size in all three
processes defined above: $\{ Z_n\}$, $\{Z^{im}_n\}$, and
$\{Z^0_n\}$. Let $g(s)$ be the immigration p.g.f.. Also, let
${\cal A}_n=\{Z_{n-1}>0\}$ for processes without immigration, and
${\cal A}_n$ be the certain event -  otherwise. The following
result is true.

\begin{theorem} If $0<m<1$, then for $x\ge 0$
\be \label{distr_limit}\lim_{n\to \infty}P(\hat{\CM}_n\le x|{\cal
A}_n)=\gamma(F(x)) \ee and \be \label{moment_limit} \lim_{n\to
\infty}E(\hat{\CM}_n|{\cal A}_n)=\sum_{k=0}^\infty
[1-\gamma(F(k))] \ee where

(i) in case of $\{Z_n\}$, $\gamma$ is the unique p.g.f. solution
of $\gamma(f(s))=m\gamma(s)+1-m$ and (\ref{moment_limit}) holds
if, in addition, $EX_i(n)\log(1+ X_i(n))<\infty$.

 (ii) in case of process $\{Z^{im}_n\}$, (\ref{distr_limit}) holds
 provided $E\log(1+Y_n)<\infty$ and
$\gamma$ is the unique p.g.f. solution of $
\gamma(s)=g(s)\gamma(f(s))$. (\ref{moment_limit}) is true if, in
addition, $EY_n<\infty$.

(iii) in case of process $\{Z^{0}_n\}$ we assume that
$E\log(1+Y_n)<\infty$. Then $ \gamma(s)= 1 - \sum_{n=0}^\infty
[1-g(f_n(s))]$ $(0< s\le 1)$ and $\gamma(0) = \{1 +
\sum_{n=0}^\infty [1-g(f_n(0))]\}^{-1}$. Also,
(\ref{moment_limit}) holds if, in addition, $EY_n<\infty$.
\end{theorem}

\begin{example} Consider $\{Z_n\}$ with geometric offspring
p.g.f. $f(s)=p/(1-qs),$\ where \ $1/2<p=1-q<1$. \ Then \
$m=q/p<1$\ and it is not difficult to see that
$\gamma(s)=(1-m)s/(1-ms).$\ Hence
\[
\lim_{n\to\i}P(\CM_n\leq k\mid Z_{n-1}>0)= \frac{\D
(p-q)(1-q^{k+1})}{\D p-q(1-q^{k+1})}.
\]
It can also be seen (\cite{RahYan99}) that
\[
\frac{m}{1-pm}\le \lim_{n\to \infty}E(\CM_n|Z_n>0)\le
\frac{m}{1-m}.
\]
\end{example}

\begin{example} Consider $\{Z^{im}_n\}$ (see \cite{Pak71}) with
\[
f(s)=(1+m-ms)^{-1} \quad (0<m<1) \quad \mbox{and} \quad  g(s)=
f^{\nu}(s) \qquad (\nu
>0).\]
Then $ \gamma(s)=((1-m)/(1-ms))^\nu$, a negative binomial p.g.f.,
and the above theorem yields
\[
  \lim_{n   \to   \infty}  P(\CM^{im}_n  \le x) =
\left( \frac{1-m}{1-mF(x)}\right)^\nu \ \mbox{and}\ \ \lim_{n\to
\infty} E\CM^{im}_n=\sum_{j=0}^{\infty} 1-\left[
\frac{1-F(j)}{1-mF(j)}\right]^{\nu}\le \frac{\nu m^2}{1-m}.
\]
\end{example}

\begin{example} Let $\mu=EY_n$. Consider $\{Z^0_n\}$ with
\[
f(s)=(1+m-ms)^{-1} \qquad \mbox{and} \qquad g(s)= 1 - (\mu /
m)\log (1+m - ms)\qquad (0<m<1).\]
  In  this case $\gamma = (m-\mu \log(1-ms))/(m-\mu \log (1-m))$ and by the
theorem
\[
\D \lim_{n   \to   \infty} P\{\CM^0_n  \le x\}  = {m  - \mu \log
(1-mF(x))  \over m -\mu \log (1-m)}\ ,
\]
and
\[
\lim_{n \to \infty}E\CM^0_n = \mu\  {\D m+\sum_{k=0}^\infty \log
\frac{1-m[(1+m)^{k+1}-m^{k+1}]}{1-m} \over m-\mu\log (1-m)}\le {\D
\mu m \over m -\mu \log(1-m)}\ {\ds m \over 1-m}.
\]
\end{example}

\subsection{Critical processes}
In the rest of this section we need some asymptotic results for
the maxima of i.i.d. random variables. Recall that a distribution
function $F(x)$ belongs to the max-domain of attraction of a
distribution function $H(x,\theta)$ (i.e., $F\in D(H)$) if and
only if there exist sequences $ a(n)>0$\ and \ $b(n)$\ such that
\begin{equation}\label{3.8}
\lim_{n\to\infty}F^n(a(n)x+b(n))=H(x, \theta)\ ,
\end{equation}
weakly. According to the classical Gnedenko's result, \
$H(x;\theta)$
 has the following (von Mises) form
\be \label{dom_attr}H(x;\theta)=\exp\{-h(x;\theta)\}
=\exp\left\{-(1+x\theta^{-1})^{-\theta}\right\}, \quad
1+x\theta^{-1}>0; \  -\infty<\theta<\infty. \ee Necessary and
sufficient conditions for \ $F\in D(H)$ are well-known. In
particular, \  $F\in D(\exp\{-x^{-a}\})$, \ $a>0$\ if and only if
for \ $x>0$ the following regularity condition on the tail
probability holds
\begin{equation} \label {24}
1-F(x)=x^{-a}L(x)\ ,
\end{equation}
where $L(x)$ is a slowly varying at infinity function (s.v.f.).

\vspace{0.5cm}{\bf A. Processes without immigration.}\ In case of
a simple BGW process, the following result holds.

\begin{theorem}  Let $m=1$ and $\sigma^2<\infty$. (i) If
(\ref{3.8}) holds, then  \be \label{thm2_dist}\lim_{n\to
\infty}P\left(\frac{\CM_n-b(n)}{a(n)}\le
x|Z_{n-1}>0\right)=\frac{1}{1+\sigma^2h(x, \theta)/2}. \ee

  (ii) If
(\ref{24}) holds, then \be  \label{thm2_exp} \lim_{n\to
\infty}\frac{E(\CM_n|Z_{n-1}>0)}{n^{1/a}L_1\left(n\right)}=
\frac{\pi/a}{\sin(\pi/a)} \qquad (a\ge 2), \ee where $L_1(x)$ is
certain s.v.f. with known asymptotics.
\end{theorem}

The theorem implies that if $F\in D(\exp\{-e^{-x}\})$ then the
limiting distribution is logistic with c.d.f.
$\left(1+e^{-x}\right)^{-1}$; and if $F\in D(\exp\{-x^{-a}\})$
then the limiting distribution is log-logistic with c.d.f.
$\left(1+x^{-a}\right)^{-1}$.

\begin{theorem}  Let $m=1$, $\sigma^2=\i$, and (\ref{24}) holds. Then for $x \geq
0$ and $1<a\le 2$ \be \label{thm3_dist}\lim_{n\to\i}P\left(
\frac{\D \CM_n}{\D n^{1/[a(a-1)]}L_2\left(n\right)} \leq x\mid
Z_{n-1}>0\right) = 1-\frac{1}{\D
\left(1+x^{a(a-1)}\right)^{1/(a-1)}} , \ee which is a Burr Type
XII distribution (e.g. \cite{Tad80}) and
\begin{equation}\label{4.77}
\lim_{n\rightarrow \infty}\frac{E(\CM_n|Z_{n-1}>0)
}{n^{1/[a(a-1)]}L_2\left(n\right)}\! \! = \frac{1}{a-1} B \left(
\frac{1}{a-1}-\frac{1}{a(a-1)}, 1+\frac{1}{a(a-1)} \right) \quad
(1<a\le 2),
\end{equation}
where $B(u,v)$ is the Beta function and $L_2(x)$ is certain s.v.f.
with known asymptotics.
\end{theorem}

Note that for $a=2$ the right-hand sides in (\ref{thm2_dist})
(under assumption (\ref{24})) and (\ref{thm2_exp}) coincide with
those in (\ref{thm3_dist}) and (\ref{4.77}), respectively. The
right-hand side in (\ref{4.77}) is the expected value of the limit
in (\ref{thm3_dist}) (see \cite{Tad80}).

\begin{example}
Let $1-F(x)\sim x^{-2}\log x$. In this case one can check (see
\cite{RahYan99}) that Theorem~3 with $a=2$ implies
\[
\lim_{n\to\i}P\left( \frac{\D \CM_n}{\D n^{1/2}(\log n)^{3/2}}\leq
x\mid Z_{n-1}>0\right) = \frac{4x^2}{1+4x^2}\ .
\]
for $x \geq 0$ and
\[
\lim_{n\rightarrow \infty}\frac{\D E(\CM_n|Z_{n-1}>0)}{\D
n^{1/2}(\log n)^{3/2}} = \frac{\pi}{2}\ .
\]
\end{example}

\vspace{0.5cm}{\bf B. Processes with immigration $\{Z^{im}_n\}$.}\
Let $\mu=EY_n$. We have the following theorem.

\begin{theorem} Assume that $ m=1, \  0<\sigma^2<\infty$, and $0 <
\mu < \infty.$ (i) If (\ref{3.8}) holds, then \be
\label{thm4_dist} \lim_{n\to \infty}P\left(
\frac{\CM^{im}_n-b(n)}{a(n)} \leq x \right) = \frac{1}{\D
(1+\sigma^2 h(x, \theta)/2)^{2\mu/\sigma^2}}.\ee (ii) If
(\ref{24}) is true, then \be \label{thm4_exp} \lim_{n\to
\infty}\frac{\D EM^{im}_n}{\D n^{1/a}L_2(n)} =
\frac{2\mu}{\sigma^2} B\left(\frac{2\mu}{\sigma^2}+\frac{1}{a},
1-\frac{1}{a}\right) \quad (a\ge 2),\ee where $B(u,v)$ is the Beta
function and $L_2(x)$ is certain s.v.f. with known asymptotics.
\end{theorem}

The theorem implies that if $F\in D(\exp\{-e^{-x}\})$ then the
limiting distribution is generalized logistic  with c.d.f.
$\left(1+\sigma^2 e^{-x}/2\right)^{-2\mu/\sigma^2}$; if $F\in
D(\exp\{-x^{-a}\})$ then the limiting distribution is a Burr Type
III (e.g. \cite{Tad80}) with c.d.f.
$\left(1+\sigma^2x^{-a}/2\right)^{-2\mu/\sigma^2}$. The right-hand
side in (\ref{thm4_exp}) is the expected value of the limit in
(\ref{thm4_dist}) (see \cite{Tad80}).

\begin{theorem}  Let $m=1$, $\sigma^2=\i$, and (\ref{24}) holds.
In addition, suppose \be \label{theta_cond}\Theta(x):=-\int_0^x
\log[1-P(Z^{im}_t>0)]dt=c\log x+d+\varepsilon(x), \ee where
$\lim_{x\to \infty}\varepsilon(x)=0$, $c>0$, and $d$ are
constants. Then for $x \geq 0$, \be
\label{thm5_dist}\lim_{n\to\i}P\left( \frac{\D \CM^{im}_n}{\D
n^{1/[a(a-1)]}L_2\left(n\right)} \leq x\right) =
\frac{1}{(1+x^{-a(a-1)})^c} \quad (1<a\le 2), \ee which is a Burr
Type III distribution (e.g. \cite{Tad80}) and \be
\label{thm5_exp}\lim_{n\rightarrow \infty}\frac{E\CM_n^{im}
}{n^{1/[a(a-1)]}L_2\left(n\right)}= cB\left(c+\frac{1}{a(a-1)},
1-\frac{1}{a(a-1)}\right)\quad (1<a\le 2), \ee where $B(u,v)$ is
the Beta function and $L_2(x)$ is certain s.v.f. with known
asymptotics. The right-hand side in (\ref{thm5_exp}) is the
expected value of the limit in (\ref{thm5_dist}).
\end{theorem}
Note that for $c=1$ and $a=2$ the right-hand sides in
(\ref{thm5_dist}) and (\ref{thm5_exp}) coincide with those in
(\ref{thm3_dist}) and (\ref{4.77}), respectively. The condition
(\ref{theta_cond}) holds even when the immigration mean is not
finite. Next example illustrates this point.

\begin{example} Following \cite{Pak75}, we consider offspring and
immigrants generated by
\[
f(s)=1-(1-s)(1+(a-1)(1-s))^{-1/(a-1)}\  \mbox{and}\
g(s)=\exp\{-\lambda (1-s)^{a-1}\},
\]
respectively. Then (\ref{24}) holds and (\ref{theta_cond}) yields
\[
\Theta(t)=(\lambda/(a-1))[\log t+\log (a-1)+\log (1+(a-1)t)^{-1}].
\]
Therefore,
\[
\lim_{n\to\i}P\left( \frac{\D \CM^{im}_n}{\D
n^{1/[a(a-1)]}L_3\left(n\right)} \leq x\right) =
\frac{1}{(1+x^{-a(a-1)})^{\lambda/(a-1)}} \quad (1<a\le 2),
\]
and
\[
\lim_{n\rightarrow \infty}\frac{E \CM_n^{im}
}{n^{1/[a(a-1)]}L_3\left(n\right)}=
\frac{\lambda}{a-1}B\left(\frac{\lambda }{a-1}-\frac{1}{a(a-1)},
1-\frac{1}{a(a-1)}\right) \quad (1<a\le 2),
\]
where $B(u,v)$ is the Beta function and $L_3(x)$ is certain s.v.f.
with known asymptotics.
\end{example}

 \vspace{0.5cm}{\bf C. Foster-Pakes processes
$\{Z^0_n\}$.} The following limit theorem for $\CM^0_n$ under a
non-linear normalization holds.

\begin{theorem} Assume that $ m=1, \  0<\sigma^2<\infty$, and $0 <
\mu < \infty.$ If \be\label{cond} \lim_{n\to\infty} {\ds
P(X_1(1)>n) \over P(X_1(1)>n+1)} = 1\  \ee then for \ $0<x <1$,
\be\label{crthm} \lim_{n\to\infty} P\left( {\D \log U(\CM^0_n)
\over \log n} \leq x \right) = x, \ee where $U(y)=1/(1-F(y))$.
\end{theorem}

Note that (\ref{cond}) is a necessary condition for $X_1(n)$ to be
in a max-domain of attraction.

\subsection{Supercritical processes}
Denote by $\hat{\CM}_n$ (as in the subcritical case above) the
maximum family size in all three processes: $\{ Z_n\}$,
$\{Z^{im}_n\}$, and $\{Z^0_n\}$. The following result is true.

\begin{theorem} Assume that $m>1$ and $EX_i(n)\log(1+ X_i(n))<\infty$. If (\ref{3.8}) holds, then
\[
\lim_{n\to \infty}P\left( \frac{\hat{\CM}_n-b(m^n)}{a(m^n)} \leq x
\right) =\psi(h(x, \theta))\] If (\ref{24}) is true, then

\[\lim_{n\to
\infty}\frac{E\hat{\CM}_n}{m^{-n/a}L_1\left(m^{-n/a}\right)}=\int_0^\infty
1-\psi(x^{-a})dx,
\]
where $L_1(x)$ is certain s.v.f. with known asymptotics.

(i) in case of $\{Z_n\}$, $\psi$ is the unique, among the Laplace
transforms, solution of \be \label{psi_eqn}
\psi(u)=f(\psi(um^{-1})), \qquad (u>0).\ee

(ii) in case of $\{Z^{im}_n\}$, we assume in addition that $E
\log(1+ Y_n) < \infty$ and \[
 \psi(u)   = \prod_{k=1}^{\infty}g(\varphi(um^{-k})) \qquad (u>0) \ ,
\] where $\varphi(u)$ is the unique, among the
Laplace transforms, solution of (\ref{psi_eqn}).

(iii) in case of $\{Z^0_n\}$,  we assume in addition that $E Y_n <
\infty$ and \[ \psi(u) = g(\varphi(u))  -  \sum_{n=0}^\infty [1 -
f(\varphi({u m^{-n}}))]P(Z^0_n=0) \qquad (u>0)\] and $\varphi(u)$
is the unique, among the Laplace transforms, solution of
(\ref{psi_eqn}). \end{theorem}

It is interesting to compare the limiting behavior of the maximum
family size in the processes allowing immigration with that when
the processes evolve in "isolation", i.e., without immigration. In
the supercritical case, as might be expected, the immigration has
little effect on the asymptotics of the maximum family size. The
limits differ only in the form of the Laplace transform $\psi(u)$.
In the subcritical and critical cases the mechanism of immigration
eliminates the conditioning on non--extinction. Theorem~6 for the
Foster-Pakes process differs from the rest of the results by the
non-linear norming of $\CM_n$. The study of the limiting behavior
of the expectation in this case needs additional efforts.

It is known that some of the most popular discrete distributions,
like geometric and Poisson, do not belong to any max-domain of
attraction. This restricts the applicability of the results in the
critical and supercritical cases above. A general construction of
discrete distributions attracted in a max-domain is given in Wilms
(1994). As it is proved there, if $X$ is attracted by a Gumbel or
Fr\'{e}chet distributions, then the same holds for the integer
part $[X]$. Next we follow a different approach considering
triangular arrays of geometric variables which leads to branching
processes with varying environments.

The results in this section are published in \cite{Mit98},
\cite{MitYan99}, and \cite{RahYan96}-\cite{RahYan99}. In
\cite{YanTso00} an extension for order statistics is considered.

\section{Maximum family size in processes with
varying environments}

It is well-known that the geometric law is not attracted to any
max-stable law. Therefore, the limit theorems for maxima in the
critical and supercritical cases above do not apply to geometric
offspring. In this section we utilize a triangular array of
zero-modified geometric (ZMG) offspring distributions, instead.

\subsection{Maxima of arrays of zero-modified geometric variables}
In this subsection we prove limit theorems for maximum of ZMG with
p.m.f. \nbeq P(X_i(n)=j)=\cases{a_np_n(1-p_n)^{j-1} & if $j\ge 1$,
\cr 1-a_n & if $j=0$,   $\qquad (n=1,2,\ldots)$} \neeq For a
 positive integer $\nu_n$ consider the triangular array of variables
\nbeq
 X_1(1), X_2(1), & \ldots, & X_{\nu_1}(1)\\
X_1(2), X_2(2), & \ldots, & \qquad X_{\nu_2}(2) \\
& \ldots & \\
X_1(n), X_2(n), &  \ldots, & \qquad \qquad \qquad X_{\nu_n}(n)
\neeq We prove limit theorems as $\nu_n\to \infty$ for the row
maxima
\[
\CM_n=\max_{1\le i\le \nu_n}X_i(n).
\]
Let $\Lambda$ has the standard Gumbel  law with c.d.f. $
\exp(-e^{-x})$ for $-\infty<x<\infty$.

\begin{theorem} Assume that for some real
$c$
\[
\lim_{n\to \infty}p_n=0 \quad \mbox{and} \quad \lim_{n\to
\infty}p_n\log(\nu_na_n)=2c.
\]

A. If $\lim_{n\to \infty}\log(\nu_na_n)=\infty$, then $c\ge0$ and
\[
p_n\CM_n-\log(\nu_na_n)\st{d}\to \Lambda -c.
\]

B. If $\lim_{n\to \infty}\log(\nu_na_n)=\alpha$,
$(-\infty<\alpha<\infty)$, then
\[
p_n\CM_n \st{d}\to (\Lambda+\alpha)^+.
\]
\end{theorem}

The idea of the proof is to exploit: (i) the exponential
approximation to the zero-modified geometric law when its mean
$a_n/p_n$ is large; (ii) the fact that exponential law is
attracted by Gumbel distribution.

\subsection{Processes with varying geometric environments}
Consider a branching process with ZMG offspring law defined over
the triangular array above. Thus, we have a simple branching
process with geometric varying environments. For this process we
prove limit theorems for the offspring maxima in all three
classes: subcritical, critical, and supercritical. Define
$\mu_0=1$,
\[
\mu_n=E(Z_n|Z_0=1)=\prod_{j=1}^nm_j \qquad (n\ge1).
\]
If the environments are weakly varying, i.e., $\mu=\lim_{n\to
\infty}\mu_n$ exists, then the processes can be classify (see
\cite{MitPakYan03}) as follows.

\nbeq \vspace{1cm} \{Z_n\}\  \mbox{is}\ \
\cases{\mbox{supercritical} & if $\mu=\infty$ \quad \quad \ \
\mbox{i.e.} \ $\sum_n(m_n-1)\to \infty$ \cr \mbox{critical} & if
$\mu\in(0,\infty)$ \quad  \mbox{i.e.} \ $\sum_n(m_n-1)< \infty$
\cr \mbox{subcritical} & if $\mu=0$ \quad \quad \quad \
\mbox{i.e.} \ $\sum_n(m_n-1)\to -\infty$} \neeq Define the maximum
family size for the process with varying geometric environments as
\[
\CM_n^{ge}=\max_{1\le i\le Z_n}X_i(n), \qquad (n=1,2,\ldots)
\]
In the result below the role played by $\nu_n$ before is played by
$B_{n-1}$ where
\[
B_n=\mu_n\sum_{j=1}^n \frac{p_j^{-1}-1}{\mu_j}.
\]
Let $\CV$ be a standard logistic random variable with c.d.f.
$(1+e^{-x})^{-1}$ for $-\infty<x<\infty$.

\begin{theorem} Suppose that $\lim_{n\to
\infty} B_n=\infty$ and for $c$ real
\[
\lim_{n\to \infty}p_n=0 \quad \mbox{and} \quad \lim_{n\to
\infty}p_n\log(B_{n-1}a_n)=2c.
\]

A. If $\lim_{n\to \infty}\log(B_{n-1}a_n)=\infty$, then
\[
(p_n\CM_n^{ge}-\log(B_{n-1}a_n)|Z_{n-1}>0)\st{d}\to \CV -c.
\]

B. If $\lim_{n\to \infty}\log(B_{n-1}a_n)=\alpha$,
$(-\infty<\alpha<\infty)$, then
\[
(p_n\CM_n^{ge}|Z_{n-1}>0) \st{d}\to (\CV+\alpha)^+.
\]
\end{theorem}

Referring to the above theorem, we can say that the branching
mechanism transforms Gumbel to logistic distribution. It is
interesting to notice that this is in parallel with results for
maximum of i.i.d. random variables with random geometrically
distributed index discussed in \cite{GneGne82}.

\begin{example} Let us sample a
linear birth and death process $(\CB_t)$ at irregular times. Let
$Z_n=\CB_{t_n}$ where $0<t_n<t_{n+1}\to t_{\infty}\le \infty$. If
$\lambda$ and $\mu$ are the birth and death rates, respectively,
and $d_n=t_n-t_{n-1}$, then $a_n=m_np_n$,
$$p_n=\cases{{\D \l-\mu\over \D \lambda m_n-\mu} & if $\lambda \not=\mu$,\cr
\frac{\D 1}{\D 1+\lambda d_n} & if $\l=\mu$,} \quad
m_n=e^{(\lambda -\mu)d_n}.$$ and
$$B_n=\cases{{\D \lambda (\mu_n-1)\over \D \lambda -\mu} & if $\lambda \not=\mu$,\cr
\lambda t_n & if $\lambda =\mu$,} \quad
\mu_n=e^{(\lambda-\mu)t_n}.$$

A. If $\lambda>\mu$ and
$$\lim_{n\to
\infty}\frac{\D t_n}{\D m_n}=\frac{2c}{\lambda-\mu}\in
[0,\infty),$$ then
\[
\left(\frac{\CM_n^{ge}}{m_n}-(\lambda-\mu)t_n)\ |\
Z_{n-1}>0\right)\st{d}\to \CV -c.
\]
B. If $\lambda=\mu$ and $t_n=n^{\delta}l(n)$ \ \ $(\delta \ge 1)$,
then
\[
\left(\frac{\CM_n^{ge}}{\lambda\delta n^{\delta-1}l(n)}-\log n\ |\
Z_{n-1}>0\right) \st{d}\to \CV.
\]
\end{example}

The results in this section can be found in \cite{MitPakYan03}.

\section{Maxima in bisexual processes}

In this section we consider maxima of triangular arrays of
bivariate geometric random vectors. The obtained results are
applied to a class of bisexual branching processes.

\subsection{Max-domain of attraction of bivariate geometric arrays}
The following construction is due to Marshall and Olkin
\cite{MO85}. Consider a random vector $(U, V)$ having Bernoulli
marginals, i.e., it takes on four possible values (0,0), (0,1),
(1,0), and (1,1) with probabilities $p_{00}, \ p_{01}, \ p_{10}$,
and $p_{11}$, respectively. Thus the marginal probabilities for
$U$ and $V$ are \nbeq P(U=0)=p_{0+}=p_{00}+p_{01}, & &
P(U=1)=p_{1+}=p_{10}+p_{11} \\
    P(V=0)=p_{+0}=p_{00}+p_{10}, & & P(V=1)=p_{+1}=p_{01}+p_{11} .
    \neeq
Consider a sequence $\{(U_n, V_n)\}_{n=1}^\infty$ of
    independent and identically distributed with $(U, V)$ random
    vectors. Let $\xi$ and $\eta$ be the number of zeros preceding the
    first 1 in the sequences $\{U_n\}_{n=1}^\infty$ and
    $\{V_m\}_{n=1}^\infty$, respectively. Both $\xi$ and $\eta$ follow
    a geometric distribution and, in general, they are dependent variables.
    The vector $(\xi, \eta)$ has a bivariate geometric distribution
    with probability mass function for integer $l$ and $k$

\begin{equation}\label{def1}
P(\xi=l, \eta=k) =
    \left\{
\begin{array}{ll}
    p_{00}^lp_{10}p_{+0}^{k-l-1}p_{+1} & \mbox{if} \quad 0\leq l <k,\\
    p_{00}^l p_{11} & \mbox{if} \quad l=k, \\
    p_{00}^k p_{01} p_{0+}^{l-k-1}p_{1+} & \mbox{if} \quad 0\leq k <l.\\
\end{array}
\right.
\end{equation}
and
\begin{equation}\label{def2}
P(\xi > l, \eta > k)=
 \left\{
\begin{array}{ll}
    p_{00}^{l+1}p_{+0}^{k-l}& \mbox{if} \quad 0\leq l \leq k,\\
    p_{00}^{k+1} p_{0+}^{l-k} & \mbox{if} \quad 0\leq k <l.\\
\end{array}
\right.
\end{equation}
The marginals of $\xi$ and $\eta$ for integer $l$ and $k$ are $
P(\xi=l)=p_{1+}p_{0+}^l \quad (l\geq 0)$ and
$P(\eta=k)=p_{+1}p_{+0}^k \quad (k\geq 0)$, respectively  and \beq
\label{marg} \bar{F}_\xi(l)=P(\xi > l)=p_{0+}^{l+1} \quad (l\geq
0), \qquad \bar{F}_\eta(k)=P(\eta > k)=p_{+0}^{k+1} \quad (k\geq
0). \eeq

For $n=1,2,\ldots$, let $\nu_n$ be a positive integer and
$\{(\xi_i(n), \eta_i(n)): i=1,2,\ldots, \nu_n\}$ be a triangular
array of independent random vectors with the same bivariate
geometric distribution (\ref{def1}) where $p_{ij}$ are replaced by
$p_{ij}(n)\ \ (i,j=0,1)$ for $n=1,2,\ldots$ That is, \nbeq
(\xi_1(1), \eta_1(1)), (\xi_2(1), \eta_2(1)), & \ldots, &
(\xi_{\nu_1}(1), \eta_{\nu_1}(1))\\
(\xi_1(2), \eta_1(2)), (\xi_2(2), \eta_2(2)), & \ldots, & \qquad
(\xi_{\nu_2}(2), \eta_{\nu_2}(2)) \\
& \ldots & \\
(\xi_1(n), \eta_1(n)), (\xi_2(n), \eta_2(n)), & \ldots, & \qquad
\qquad \qquad  (\xi_{\nu_n}(n), \eta_{\nu_n}(n)) \neeq Below we
prove a limit theorem as $\nu_n \to \infty$ for the bivariate row
maximum
\[
(\CM_n^\xi, \CM_n^\eta)=\left(\max_{1\leq i\leq \nu_n} \xi_i(n),
\max_{1\leq i\leq \nu_n}\eta_i(n)\right).
\]

\begin{theorem} Let $\lim_{n\to \infty}\nu_n=\infty$. If there are
constants $0\leq a, b, c < \infty$, such that \be \label{assum2}
\lim_{n\to \infty}p_{11}(n)\log \nu_n = 2c \quad \lim_{n\to
\infty} \frac{p_{10}(n)}{p_{11}(n)}\log \nu_n = a \quad \mbox{and}
\quad \lim_{n\to \infty} \frac{p_{01}(n)}{p_{11}(n)}\log \nu_n =
b, \ee then for
$x, y \geq 0$ 
\nbeq \lefteqn{ \lim_{n\to\infty}
    P\left( p_{11}(n)\CM_n^\xi-\log\nu_n\leq x, \ \ p_{11}(n)\CM_n^\eta-\log\nu_n\leq
      y \right)}\\
      & &  =
    \exp\left\{ -e^{ -x-a-c}-e^{ -y-b-c}+e^{ -\max\{x,y\}-a -b -c}\right\}.
    \neeq
\end{theorem}

{\bf Proof}\ Set $x_n=(x+\log \nu_n)/p_{11}(n)$ and $y_n=(y+\log
\nu_n)/p_{11}(n)$.  \nbeq P\left( \CM_n^\xi\leq x_n,
\CM_n^\eta\leq
y_n\right) & = & (F(x_n, y_n))^{\nu_n} \\
    & = &
    \left(1-\bar{F}_\xi(x_n)-\bar{F}_\eta(y_n)+P(\xi_i(n)>x_n, \eta_i(n)>y_n)\right)^{\nu_n}.
    \neeq
   Let $x<y$ and thus, $x_n<y_n$. Taking logarithm, expanding in Taylor
   series,
   and using
    (\ref{def2}) and (\ref{marg}), we obtain
 \beq \label{lim0}
 \lefteqn{\log P\left( \CM_n^\xi\leq x_n, \CM_n^\eta\leq
y_n\right)}\\
 & = & \nu_n \log
\left(1-\bar{F}_\xi(x_n)-\bar{F}_\eta(y_n)+P(\xi_i(n)>x_n, \eta_i(n)>y_n)\right) \nonumber \\
    & = &
    -\nu_n \left\{ [\bar{F}_\xi(x_n)+\bar{F}_\eta(y_n)-P(\xi_i(n)>x_n, \eta_i(n)>y_n)](1+o(1))\right\}\nonumber \\
    & = &
    -\left( \nu_n p_{0+}(n)^{[x_n]+1} + \nu_n p_{+0}(n)^{[y_n]+1}
    - \nu_n p_{00}(n)^{[x_n]+1}p_{+0}(n)^{[y_n]-[x_n]}\right)
    (1+o(1)) . \nonumber
    \eeq
    Write $[x_n]=x_n-\{x_n\}$, where $0\leq \{x_n\}<1$ is
    the fractional part of $x_n$. It is easily seen that
    $\lim_{n\to \infty}(p_{0+}(n))^{[x_n]+1}=\lim_{n\to
    \infty}(p_{0+}(n))^{x_n+1-\{x_n\}}=
    \lim_{n\to \infty}(p_{0+}(n))^{x_n}$ as $n\to \infty$.
    Furthermore, taking into account (\ref{assum2}), we have
    \nbeq
    \lefteqn{\log \left(\nu_n p_{0+}^{x_n}(n)\right)=
    \log \nu_n + \frac{x+\log
    \nu_n}{p_{11}(n)}\log(1-p_{1+}(n))}\nonumber \\
    & = &
    \log \nu_n - \frac{x+\log
    \nu_n}{p_{11}(n)}\left(p_{11}(n)+p_{10}(n)
    +\frac{1}{2}(p_{11}(n)+p_{10}(n))^2+O(p^3_{1+}(n))\right)
        \nonumber \\
    & = &
    -x(1+o(1)) - \frac{p_{10}(n)}{p_{11}(n)}\log \nu_n -
    \frac{(p_{11}(n)+p_{10}(n))^2}{2p_{11}(n)}\log \nu_n +
    O(p^2_{11}(n))
    \nonumber \\
    & = &
    -x(1+o(1))-\left(\frac{p_{10}(n)}{p_{11}(n)}+
    \frac{1}{2}p_{11}(n)\right)\log\nu_n (1+o(1)) +
    O(p^2_{11}(n))
    \nonumber \\
    & \to &
    -x -a -c \ .
    \nonumber
    \neeq
    Therefore
    \be \label{lim1}
    \lim_{n \to \infty}\nu_n p_{0+}(n)^{[x_n]+1}=e^{\D -x -a -c} \ .
    \ee
    Similarly we arrive at
    \be \label{lim2}
    \hspace{-0.3cm} \lim_{n \to \infty}\nu_n p_{+0}(n)^{[y_n]+1}=e^{\D -y -b -c} \
    \mbox{and} \
    \lim_{n \to \infty}\nu_n p_{00}(n)^{[x_n]+1}=e^{\D -x -a -b -c} \ .
    \ee
    Finally,
    \nbeq
    \lefteqn{\log p_{+0}^{y_n-x_n}(n)=
    \frac{(y-\log\nu_n)-(x-\log\nu_n)}{p_{11}(n)}\log
    (1-p_{11}(n)-p_{01}(n))} \\
    & = &
    -\frac{y-x}{p_{11}(n)}\left(p_{11}(n)+p_{01}(n)
    +\frac{1}{2}(p_{11}(n)+p_{01}(n))^2+O(p^3_{+1}(n))\right)\\
    & = &
    x-y -(y-x)\left(
    \frac{p_{01}(n)}{p_{11}(n)}(1+o(1))+\frac{1}{2}p_{11}(n)(1+o(1))
    + O(p^2_{+1}(n))\right)\\
    & \to &
    x-y
    \neeq
    Thus,
    \be \label{lim3}
    \lim_{n \to \infty}p_{+0}(n)^{[y_n]-[x_n]}=e^{\D x-y} \ .
    \ee
    The assertion of the theorem for $x<y$ follows from
    (\ref{lim0})-(\ref{lim3}). The case $y<x$ is treated
    similarly. This completes the proof.

    In particular, if $a=b=0$ then
    \[
\lim_{n\to\infty}
    P\left( p_{11}(n)\CM_n^\xi-\log\nu_n\leq x, p_{11}(n)\CM_n^\eta-\log\nu_n\leq
      y \right)
     =
    \exp\left\{ -e^{ -\min\{x,y\}-c}\right\}.
    \]
    Note that in this case the limit is proportional to the upper bound  for the possible asymptotic
    distribution of a multivariate maximum given in
    \cite{Gal87}, Theorem~5.4.1.

   For the componentwise maxima, applying Theorem~10, one can obtain
   the following limiting results. If
   $p_{1+}(n)\log \nu_n \to 2c_1<\infty$, then
    \[
    \lim_{n\to\infty}P( p_{11}(n)\CM_n^\xi-\log\nu_n\leq x)
 =
    \exp\left\{ -e^{ -x-c_1}\right\}.
    \]
    If $p_{1+}(n)\log \nu_n \to 2c_2<\infty$, then
    \[
    \lim_{n\to\infty}
    P\left(p_{11}(n)\CM_n^\eta-\log\nu_n\leq
      y \right) =
    \exp\left\{ -e^{ -y-c_2}\right\}.
    \]

\subsection{Bisexual processes with varying geometric
environments}  Consider the array of bivariate random vectors
$\{(\xi_i(n),\eta_i(n)):\ i=1,2,\ldots; \ n=0,1,\ldots\}$, which
are independent with respect to both indexes. Let
$L:\CR^+\times\CR^+\to\CR^+$ be a mating function. A bisexual
process with varying environments is defined (see \cite{MMR04}) by
the recurrence: $Z_0=N>0$,
\[
(Z^F_{n+1},Z^M_{n+1})=\sum_{i=1}^{Z_n}(\xi_i(n),\eta_i(n))\] and
\[
Z_{n+1}=L(Z^F_{n+1},Z^M_{n+1}) \quad (n=0,1,\ldots).
\]
Define the mean growth rate per mating unit
\[
r_{nj}=j^{-1}E(Z_{n+1}|Z_n=j) \quad (j=1,2,\ldots) \quad
\mbox{and} \quad \mu_n=\prod_{i=0}^{n-1}r_{i1}, \ \mu_0=1 \quad
(n=1,2, \ldots)
\]

{\bf Lemma } (\cite{MMR04}) {\it If
 \be \label{assump1}
 \sum_{n=0}^\infty \left( 1- \frac{r_{n1}}{r_n}\right)\
 \ee
 then
\[
 \lim_{n\to \infty} \frac{Z_n}{\mu_n
 } = W \quad \mbox{a.s.},
 \]
 where $W$ is a nonnegative random variable with
 $E(W)<\infty$.

 If, in addition, there exist constants $A>0$ and $c>1$ such that
 \be \label{assump2}
 \prod_{i=j}^{n+j-1}r_{i1}\geq Ac^n \quad j=1,2,\ldots; \
 n=0,1,\ldots
 \ee
 and there exists a random variable $X$ with $E(X\log(1+X))<\infty$
 such that for any $u$
 \be \label{assump3}
 P(X\leq u) \leq P\left(\frac{L(\xi_i(n), \eta_i(n))}{r_{n1}}\leq u\right)\qquad
 (n=0,1,\ldots),
 \ee
then $P(W>0)>0$.}

Further on we assume that $(\xi_i(n),\eta_i(n))$ are i.i.d. copies
of the bivariate geometric vector $(\xi, \eta)$ introduced above
and that the mating is promiscuous, i.e., \be \label{mating}
L(\xi(n),\eta(n))=\xi(n)\min\{1,\eta(n)\}.\ee

\begin{theorem} \label{cpm2} Let $\{ Z_n\}$ be a bisexual branching process with varying geometric environments
and mating function (\ref{mating}).
If  
\be \label{cpm_lim_assum1}
    \prod_{j=1}^\infty p_{+0}(j)p_{0+}(j)\neq 0 \quad \mbox{and} \quad \sum_{n= 0}^\infty
    p_{+1}(n)<\infty \ ,
    \ee
    then
\be \label{cpm_Wlim1}
 \lim_{n\to \infty} \frac{Z_n}{\mu_n
 } = W \quad \mbox{a.s.},
 \ee
 where $W$ is a nonnegative random variable with
 $E(W)<\infty$ and $P(W>0)>0$.
    \end{theorem}

{\bf Proof}\ To prove the theorem it is sufficient to verify the
assumptions (\ref{assump1})-(\ref{assump3}) in the above lemma.
First, we prove that (\ref{assump1}) holds. Indeed, for $j\geq 1$
\beq jr_{nj} & = &
E(Z^F_{n+1}\min \{1,Z^M_{n+1}\}) \label{growth rate}\\
    & = &
    EE(Z^F_{n+1}\min \{1,Z^M_{n+1}\}\ |\ Z^M_{n+1}) \nonumber \\
    & = &
    (1-P(Z^M_{n+1}=0))EZ^F_{n+1} \nonumber \\
    & = &
    (1-p_{+1}^j(n))\frac{jp_{0+}(n)}{p_{1+}(n)}\ , \nonumber
    \eeq
    where we have used that both $Z^M_{n+1}$ and $Z^F_{n+1}$ are
    negative binomial with parameters $(j, p_{+1}(n))$
    and $(j,p_{1+}(n))$, respectively. Thus,
    \be \label{cpm_rn}
    r_n=\lim_{j\to \infty}r_{nj}=\lim_{j\to
    \infty}
    (1-p_{+1}^j(n))\frac{p_{0+}(n)}{p_{1+}(n)}=\frac{p_{0+}(n)}{p_{1+}(n)}
    \ .
    \ee
    Now, (\ref{growth rate}) and (\ref{cpm_rn}) imply
    $
    1-r_{n1}/r_n=p_{+1}(n)
    $,
    which along with (\ref{cpm_lim_assum1}) leads to
    (\ref{assump1}).

    Let us prove (\ref{assump3}). Indeed, for $k\geq 1$
    \nbeq
    P(L(\xi(n),\eta(n))=k)\! & \!= \!& \! \sum_{j=1}^\infty
    P(\xi(n)\min\{1,\eta(n)\}=k|\eta(n)=j)P(\eta(n)=j) \\
     \!   & = & \! P(\xi(n)=k)\sum_{j=1}^\infty P(\eta(n)=j) \nonumber \\
   \! & = & \! p_{1+}(n)p_{0+}^k(n)\sum_{j=1}^\infty
    p_{+1}(n)p_{+0}^j(n) \nonumber \\
  \!  & = & \!
    p_{+0}(n)p_{1+}(n)p_{0+}^k(n)\ .
        \nonumber
    \neeq
Therefore,
    $
    P(L(\xi(n),\eta(n))/r_{n1}\geq u)=p_{0+}^{[ur_{n1}]+1}(n)
    $
    and hence, similarly to (\ref{lim1}), taking into account
    (\ref{cpm_lim_assum1}), we obtain
    \nbeq
    \log P\left(\frac{\D L(\xi(n),\eta(n))}{\D r_{n1}}\geq u\right) & \sim & ur_{n1}\log p_{0+}(n) \\
        & = &
        -u
        \frac{p_{+0}(n)p_{0+}(n)}{p_{1+}(n)}p_{1+}(n)(1+o(1))
        \\
        & \to & -u
        \neeq
        Thus,
    $\lim_{n\to \infty}P(L(\xi(n),\eta(n)/r_{n1})\geq u)=e^{-u}$, which implies (\ref{assump3}).

    Finally, to prove (\ref{assump2}), observe that (\ref{growth rate}) implies for any $j$
    and $n$
    \nbeq
    \prod_{i=j}^{n+j-1}r_{i1} & = &
    \prod_{i=j}^{n+j-1}p_{+0}(i)p_{0+}(i)
    \prod_{i=j}^{n+j-1}p_{11}^{-1}(i) \\
        & \geq &
    \prod_{i=1}^\infty p_{+0}(i)p_{0+}(i)
    \prod_{i=j}^{n+j-1}p_{11}^{-1}(i) \\
    & \geq &
        Ac^n \ ,
    \neeq
    where $A=\prod_{i=1}^\infty p_{+0}(i)p_{0+}(i)>0$ (provided that the product in (\ref{cpm_lim_assum1}) is finite) and
    $c=\min_{i\geq j}p^{-1}_{11}(i)>1$ ($p_{11}(i)\to 0$ under
    (\ref{cpm_lim_assum1})). (\ref{assump2}) also holds if
    the product in (\ref{cpm_lim_assum1}) is infinite. Now, referring to the above lemma we complete the
    proof of the theorem.

Define offspring maxima in the bisexual process $\{Z_n\}$ by
\[
(\CM_n^F, \CM_n^M)=\left(\max_{1\le i\le Z_{n}}\xi_i(n),\
\max_{1\le i\le Z_{n}}\eta_i(n)\right).
\]

\begin{theorem} Assume that $\mu_n\to \infty$ and
there are constants $0\leq a, b, c < \infty$, such that \be
\label{assum_1} \lim_{n\to \infty}p_{11}(n)\log \mu_n = 2c \quad
\lim_{n\to \infty} \frac{p_{10}(n)}{p_{11}(n)}\log \mu_n = a \quad
\mbox{and} \quad \lim_{n\to \infty}
\frac{p_{01}(n)}{p_{11}(n)}\log \mu_n = b. \ee Also assume that
\be \label{assum_2}
    \prod_{j=1}^\infty p_{+0}(j)p_{0+}(j)\neq 0 \quad \mbox{and} \quad \sum_{n= 0}^\infty
    p_{+1}(n)<\infty \ .
    \ee
Then
\[
\lim_{n\to\infty}
    P\left( p_{11}(n)\CM_n^F -\log\mu_n\le x, p_{11}(n)\CM_n^M-\log\mu_n \le y\right)
= \int_0^\infty (G(x,y))^zdP(W\le z),
\]
where \[ G(x,y)=\exp\left\{ -e^{ -x-a-c}-e^{ -y-b-c}+e^{
-\max\{x,y\}-a -b -c}\right\}.
    \]
\end{theorem}

{\bf Proof}\ Set $x_n=(x+\log \mu_n)/p_{11}(n)$ and $y_n=(y+\log
\mu_n)/p_{11}(n)$. Under assumption (\ref{assum_1}), Theorem 11
implies \be \label{theorem_11} P\left( \CM_n^F\leq x_n,
\CM_n^M\leq y_n \ | \ Z_n=k\right) = (F(x_n, y_n))^{k}\to H(x,y).
\ee Under (\ref{assum_2}), Theorem 12 implies \be
\label{theorem_12} \lim_{n\to \infty}P\left( \frac{Z_n}{\mu_n}\le
x\right)=P(W\le x). \ee Therefore, by (\ref{theorem_11}) and
(\ref{theorem_12}),
\begin{eqnarray*}
\lefteqn{\hspace{-2cm}P\left(\CM_n^F\leq \frac{x+\log
\mu_n}{p_{11}(n)}, \CM_n^M\le \frac{y+\log
\mu_n}{p_{11}(n)}\right)=\sum_{k= 0}^\infty P\left(Z_n=k
 \right)(F(x_n,
y_n))^{k}}\\
& = & \sum_{k=0}^\infty P\left(
\frac{Z_n}{\mu_n}=\frac{k}{\mu_n}\right)(F(x_n, y_n))^{\mu_n
k/\mu_n} \\
& = & \int_0^\infty (G(x,y))^zdP(W\le z).
\end{eqnarray*}

Next example, adopted from \cite{Mit05}, shows that the various
conditions in Theorem~13 can be satisfied.
\begin{example} Let $\alpha>1$ and $\beta>1$. Set
\[
p_{11}(n)=n^{-\alpha} \qquad \mbox{and}\qquad
p_{01}(n)=p_{10}(n)=n^{-(\alpha+\beta)} \qquad (n\ge 2). \] It is
not difficult to see that with this choice of $p_{ij}(n)$
$(i,j=0,1)$, we have
\[
\log \mu_n\sim \alpha n \log n \qquad \mbox{as}\qquad n\to \infty
\]
and both (\ref{assum_1}) (with $a=b=c=0$) and (\ref{assum_2}) are
satisfied.
\end{example}

The exposition in this section follows \cite{Mit05}, extending
some of the results there.

\section{Maximum score}

In this section we assume that every individual in a Galton-Watson
family tree has a continuous random characteristic which maximum
is of interest.

\subsection{Maximum scores in Galton-Watson processes}

Let us go back to the simple BGW process and attach random scores
to each individual in the family tree. More specifically,
associate with the $j$-th individual in the $n$-th generation a
continuous random variable $Y_j(n)$. Arnold and Villase\~{n}or
(1996) published the first paper studying the maxima individual
scores ("heights"). Pakes (1998) proves more general results
concerning the laws of offspring score order statistics. Quoting
\cite{Pak98}, "these results provide examples of the behavior of
extreme order statistics of observations from samples of random
size." Define by $M_{(k),n}$ the $k$-th largest score within the
$n$-th generation and by $\bar{M}_{(k),n}$ the $k$-th largest
among the random variables $\{Y_{i}(n):\ 1\le i\le Z_\nu, 0\le
\nu\le n\}$, i.e., the $k$-th largest score up to and including
the $n$-th generation. Pakes (1998) studies the limiting behavior
of "near maxima", i.e., (upper) extreme order statistics
$M_{(k),n}$ and $\bar{M}_{(k),n}$ when $n\to \infty$ and $k$
remains fixed. The two general cases that arise are whether the
law of $Z_n$ (or the total progeny $T_n=\sum_{\nu=0}^nZ_\nu$),
conditional on survival, do not require or do require,
normalization to converge to non-degenerate limits.

If no normalization is required then no particular restriction
need to be placed on the score distribution function $S$, but the
limit laws are rather complex mixtures of the laws of extreme
order statistics. The principal result states that
\[
\lim_{n\to\infty}P(M_{(k),n}\le x| {\cal
A}_n)=\sum_{j=1}^\infty\sum_{i=0}^{k-1}{j \choose
i}(1-S(x))^iS^{j-i}(x)g_j,
\]
where it is assumed that the conditional law ${\cal G}_n$ of $Z_n$
given ${\cal A}_n$ (${\cal A}_n$ includes non-extinction)
converges to a discrete and non-defective limit ${\cal G}$ and
$g_j$ denote the masses attributed to $j$ by ${\cal G}$.

If normalization is required then one must assume that the score
distribution function $S$ is attracted to an extremal law, and
then the limit laws are mixtures of the classical limiting laws of
extreme order statistics. let us assume that there are positive
constants $C_n\uparrow \infty$ such that for the conditional law
${\cal G}_n$ we have $ {\cal G}_n(xC_n)\Rightarrow N(x), $ where
$N(x)$ is a non-defective but possibly degenerate distribution
function. Assume also that the score distribution function $S$ is
in the domain of attraction of on extremal law given by
(\ref{dom_attr}). The general result in \cite{Pak98} is \be
\label{norm_limit}\lim_{n\to\infty}P\left(\frac{M_{(k),n}-b(C_n)}{a(C_n)}\le
x|{\cal A}_n\right)=\sum_{i=0}^{k-1}\frac{(h(x,
\theta))^i}{i!}\int_0^\infty y^ie^{-yh(x, \theta)}dN(y). \ee

\begin{example} Consider an immortal (i.e., $P(X=0)=0$) supercritical process with shifted geometric offspring
law given by its p.g.f. $f(s)=s/(1+m-ms)$ $(m>1)$, then (see Pakes
(1998)) (\ref{norm_limit}) becomes
\[
\lim_{n\to\infty}P\left(\frac{M_{(k),n}-b(C_n)}{a(C_n)}\le x|{\cal
A}_n\right)=1-\left(\frac{h(x, \theta)}{1+h(x, \theta)}\right)^k.
\]
Thus the limit has a generalized logistic law when the score law
is attracted to Gumbel law, $h(x, \theta)=e^{-x}$; and a
Pareto-type law results when $S$ is attracted to the Fr\'{e}chet
law.
\end{example}

 Phatarford (see \cite{Pak98}) has raised the question (in the context
of horse racing), "What is the probability that the founder of a
family tree is better than all its descendants?" The answer turns
out to be $E\left(T^{-1}\right)$, where $T={\displaystyle
\sum_{n=0}^\infty Z_n}$ is the total number of individuals in the
family tree. More generally, if $\tau_n$ is the index of the
generation up to the $n$-th which contains the largest score,
Pakes (1998) proves that
\[
P(\tau_n=k)=E\left(\frac{Z_k}{T_n}\right), \quad (k=0,1,\ldots,n),
\]
as well as limit theorems for $\tau_n$ as $n\to \infty$.

This subsection is based on \cite{ArnVil96} and \cite{Pak98}.

\subsection{Maximum scores in two-type processes}
Let each individual in a two-type branching process be equipped
with a non-negative continuous random variable - individual score.
We present limit theorems for the maximum individual score.
Consider two independent sets of independent random vectors with
integer nonnegative components
\[\{{\bf X}^1(n)\}=\{(X^1_{1j}(n), X^1_{2j}(n))\} \ \ \mbox{and} \ \
\{{\bf X}^2(n)\}=\{(X^2_{1j}(n), X^2_{2j}(n))\} \ \  (j\ge 1; n\ge
0).
\]
A two-type branching process $\{{\bf Z}(n)\}=\{(Z_1(n),Z_2(n))\}$
is defined as follows: ${\bf Z}(0)\neq {\bf 0}$ a.s. and for
$n=1,2,\ldots$
$$ Z_1(n)=\sum_{j=1}^{Z_1(n-1)} X^1_{1j}(n) +
\sum_{j=1}^{Z_2(n-1)} X^2_{1j}(n),$$
\[Z_2(n)=\sum_{j=1}^{Z_1(n-1)} X^1_{2j}(n) + \sum_{j=1}^{Z_2(n-1)}
X^2_{2j}(n).\] Here $X^i_{kj}(n)$ refers to the number of
offspring of type $k$ produced by the $j$-th individual of type
$i$. With the $j$-th individual of type $i$ living in the $n$-th
generation we associate a non-negative continuous random variable
$\zeta_{ij}(n)$, $(i=1,2)$ "score", say. Assume that the offspring
of type $1$ and type $2$ have scores, which are independent and
identically distributed within each type. Define the maximum score
within the $n$-th generation by
\[\CM_n^\zeta =\max\{\CM_n^{\zeta_1},\ \CM_n^{\zeta_2}\}, \qquad \mbox{where}\quad
\CM_n^{\zeta_i}=\max_{1\le j\le Z_i(n-1)}\zeta_{ij}(n)\qquad
(i=1,2).
\]
Note that this is maximum of random number, independent but
non-identically distributed random variables. Let
$F_i(x)=P(\zeta_i\le x)$ $(i=1,2)$ be the c.d.f.'s of the scores
of type 1 and type 2 individuals, respectively.

\noindent{\it Assumption 1}\ (tail-equivalence)\ We assume that
$F_1$ and $F_2$ are tail equivalent, i.e., they have the same
right endpoint $x_0$ and for some $A>0$
$$ \lim_{x \uparrow x_0} \frac{1-F_1(x)}{1-F_2(x)}=A.$$

\noindent{\it Assumption 2}\ (max-stability)\ Suppose $F_1$ is in
a max-domain of attraction, i.e., (\ref{3.8}) holds.

We consider the critical branching process ${\bf Z}(n)$ with mean
matrix ${\bf M}$, which is positively regular and nonsingular. Let
${\bf M}$ has maximum eigenvalue 1 and associated right and left
eigenvectors ${\bf u}=(u_1,u_2)$ and ${\bf v}=(v_1,v_2)$,
normalized such that ${\bf u \cdot v}=1$ and ${\bf u\cdot 1}=1$.

\begin{theorem}\  Let $\{{\bf Z}(n)\}$ be the above critical
two-type branching process. If the offspring variance $2B<\infty$
and both Assumptions 1 and 2 hold, then \be \label{theorem_2D}
\hspace{-0.3cm} \lim_{n \to \infty}
P\left(\frac{\CM_n^\zeta-b(v_1Bn)}{a(v_1Bn)} \le x|{\bf Z}(n) \ne
{\bf 0} \right)=\frac{1}{1+h(x,\theta)+(v_2/v_1)h(cx+d,\theta)},
\ee where if $-\infty<\theta<\infty$ is fixed, then
$c=A^{1/|\theta |}$ and $d=0$; if $\theta\to \pm \infty$, then
$c=1$ and $d=\ln A$. \end{theorem}

 {\bf Proof}. Since $F_1(x)$ and $F_2(x)$ are
tail-equivalent, we have (see \cite{Res87}, p.67)
\[
\lim_{n\to \infty}\left(F_2(a(n)x + b(n))\right)^n \to
H(cx+d,\theta),
\]
where the constants $c$ and $d$ are as in (\ref{theorem_2D}). On
the other hand, it is well-known (see \cite{AthNey72}, p.191) that
for $x>0$ and $ y>0$ \[
 \lim_{n\to\infty}
P\left(\frac{Z_1(n)}{v_1Bn}\le x,\frac{Z_2(n)}{v_2Bn}\le y|{\bf
Z}(n) \ne {\bf 0} \right)=G(x,y), \]
 where the limiting distribution has
Laplace transform
\begin{equation}\label{LT}
   \psi(\lambda, \mu)=\frac{1}{1+\lambda+\mu}
\qquad (\lambda>0, \ \mu > 0).
\end{equation}
Set $x_n=a(v_1Bn)x+b(v_1Bn)$, $s_n=k/v_1Bn$, and $t_n=l/v_2Bn$.
Referring to the definition of both $\CM_n^\zeta$ and process
$\{{\bf Z}(n)\}$ we obtain
\begin{eqnarray*}
\lefteqn{P\left(\CM_n^\zeta \leq x_n|{\bf Z}_n\neq {\bf
0}\right)=\sum_{(k,l)={\bf 0}}^\infty P\left( {\bf
Z}(n)=(k,l)|{\bf Z}(n)\neq {\bf
 0}
 \right)
P\left(\max\left\{\CM_n^{\zeta_1}, \CM_n^{\zeta_2}\right\}\le
x_n\right)}
\\
& = & \sum_{(k,l)={\bf 0}}^\infty P\left(
\frac{Z_1(n)}{v_1Bn}=\frac{k}{v_1Bn},
\frac{Z_2(n)}{v_2Bn}=\frac{l}{v_2Bn}|{\bf Z}(n)\neq {\bf 0}
\right) \left[F_1(x_n)\right]^k\left[F_2(x_n)\right]^l
\\
&  = & \sum_{(k,l)={\bf 0}}^\infty P\left(
\frac{Z_1(n)}{v_1Bn}=s_n, \frac{Z_2(n)}{v_2Bn}=t_n|{\bf Z}(n)\neq
{\bf 0} \right)
\left[F_1(x_n)\right]^{(v_1Bn)s_n}\left[F_2(x_n)\right]^{(v_1Bn)t_n(v_2/v_1)}
\\
& \to &
\int_0^\infty \int_0^\infty H(x,\theta)^s H(cx+d,\theta)^{(v_2/v_1)t}d G(s,t) \\
& = & \int_0^\infty \int_0^\infty \exp \left\{ -sh(x,\theta)-t\frac{v_2}{v_1} h(cx+d, \theta)\right\}d G(s,t)\\
& = & \left[ 1+h(x, \theta)+\frac{v_2}{v_1} h(cx+d,
\theta)\right]^{-1},
\end{eqnarray*}
where in the last formula we used the Laplace transform of
$G(u,v)$ given in (\ref{LT}). The proof is complete.

The two examples below illustrate the kind of limit laws that can
be encountered.

\begin{example} Let $F_1$ and $F_2$ be Pareto c.d.f.'s given for \
$x_i >\theta_i>0$ and $ c>0$ by
\[
F_i(x_i)=1-\left(\frac{\theta_i}{x_i}\right)^c \qquad (i=1,2).
\]
Note that the two distributions share the same value of the
parameter $c$. It is not difficult to check that the limit is
log-logistic given by
$$\lim_{n \to \infty}
P\left\{\frac{\CM_n^\zeta}{\theta_1(v_1Bn)^{1/c}} \le x|{\bf Z}(n)
\ne {\bf 0}
\right\}=\left[1+\left(1+\frac{v_2}{v_1}\left(\frac{\theta_1}{\theta_2}\right)^{-c}\right)x^{-c}\right]^{-1}.$$
\end{example}

\begin{example} Let $F_1$ and $F_2$  be logistic and exponential
c.d.f.'s given by
\[
F_1(x_1)=1-e^{-x_1} \quad (0<x_1<\infty)\quad \mbox{and}\quad
F_2(x_2)=\frac{1}{1+e^{-x_2}} \quad (-\infty < x_2 < \infty),
\]
respectively. It is known  that both are in the max-domain of
attraction of $H(x)=\exp\{-\exp\{-x\}\}$ and share (see
\cite{AhsNev01}, p.91) the same normalizing constants $a(n)=1$ and
$b(n)=\ln n$. This fact, after inspecting the proof of the
theorem, allows us to bypass the tail-equivalence assumption and
obtain a logistic limiting distribution, i.e, for $-\infty <x
<\infty$
$$\lim_{n \to \infty} P\left\{\CM_n^\zeta-\log (v_1Bn) \le x\ |\ {\bf Z}(n)
\ne {\bf 0} \right\}=
\left[1+\left(1+\frac{v_2}{v_1}\right)e^{-x}\right]^{-1}.$$
\end{example}

The results in this subsection are modifications of those in
\cite{MitYan02}.

\section*{Acknowledgments.} I thank I. Rahimov for igniting my interest to
extremes in branching processes. Thanks to the organizers of ISCPS
2007 for the excellent conference. This work is partially
supported by NFSI-Bulgaria, MM-1101/2001.


\begin{thebibliography}{999}
\bibitem {AhsNev01} Ahsanullah, M. and Nevzorov, V. B.,
{\it Ordered Random Variables}, Nova Science Publishers, Inc.,
Huntington, NY, 2001.
\bibitem 
{ArnVil96} Arnold, B. C. and Villase\~{n}or, J. A.,
 The tallest man in the world., In: {\it Statist. Theory and
Appl.: papers in honor of Herbert A. David}, Eds. Nagaraja, H. N.,
Sen, P. K. and Morrison, D. F, pp. 81--88, Springer, Berlin, 1996.
\bibitem {AthNey72} Athreya, K. B. and Ney, P. E., {\it Branching Processes},
Springer, New York, 1972.
\bibitem {Gal87} Galambos, J. {\it The Asymptotic Theory of Extreme Order
Statistics}, 2nd Edn., Krieger, Melbourne, Florida, 1987.
\bibitem{GneGne82} Gnedenko, B.V. andGnedenko, D.B.
Laplace distributions and the logistic distribution as limit
distributions in probability theory. {\it Serdica}, 8(1982),
2:229--234 (In Russian).
\bibitem {Har63} Harris, T.E. {\it The Theory of Branching
Processes}, Springer, Berlin, 1963.
\bibitem {JagNer84} Jagers, P. and Nerman, O. Limit theorems for
sums determined by branching and exponentially growing processes.
{\it Stoch. Proc. Appl.} 17(1984), 47-71.
\bibitem{MO85}Marshall, A.W. and Olkin, I. A family of bivariate
distributions generated by the Bernoulli distribution. {\it J.
Amer. Statist. Assoc.}, 80(1985), 332-338.
\bibitem {Mit98} Mitov, K.V. The maximal number of offspring of
one particle in a branching process with state-dependent
immigration. {\it Proc. of 27th  Spring Conference  of the Union
of Bulgarian Mathematicians, Math. and Math. Education}, 1998,
92-97.
\bibitem {Mit05} Mitov, K.V. Extremes of bivariate geometic
variables with application to bisezual branching processes. {\it
Pliska Studia Math. Bulgarica}, 17(2005), 349-362.
\bibitem 
{MitPakYan03} Mitov, K.V., Pakes, A.G., and Yanev, G.P. Extremes
of geometric variables with applications to branching processes.
{\it Statist and Probab Letters}, 65(2003), 379-388.
\bibitem 
{MitYan99}  Mitov, K.V. and Yanev G.P. Maximum family size in
branching processes with state--dependent immigration, {\it Proc.
of 28th  Spring Conference  of the Union of Bulgarian
Mathematicians, Montana, Math. and Math. Education}, 1999,
142-144.
\bibitem 
{MitYan02} Mitov, K.V. and Yanev, G.P. Maximum individual score in
critical two-type branching processes. {\it C. R. Acad. Bulg.
Sci.}, 55(2002), 11:17-22.
\bibitem{MMR04} Molina, M., Mota, M., and Ramos, A. Limiting
behaviour for superadditive bisexual Galton-Watson processes in
varying environments. {\it Test}, 13(2004), 2:481-499.
\bibitem
{Pak71} Pakes, A.G. Branching processes with immigration. {\it J.
Appl. Prob.}, 8(1971), 32-42.
\bibitem{Pak75} Pakes, A.G. Some new limit theorems for the critical
branching process allowing immigration. {\it Stoch. Proc. Appl.}
3(1975), 175-185.
\bibitem 
{Pak98}  Pakes, A.G. Extreme order statistics on Galton--Watson
trees. {\it Metrika}, 47(1998), 95-117.
\bibitem {RahYan96} Rahimov, I. and Yanev G.P. On a
maximal sequence associated with simple branching processes. {\it
Institute of Mathematics and Informatics}, Sofia, Preprint no.6,
1996, pp.14.
\bibitem 
{RahYan97} Rahimov, I. and Yanev, G.P. Maximal number of direct
offspring in simple branching processes. {\it Nonlinear Analysis,
Theory, Methods and Applications}, 30(1997), 2015-2023.
\bibitem 
{RahYan99}  Rahimov, I. and Yanev, G.P. On maximum family size in
branching processes. {\it J. Appl. Probab.}, 36(1999), 632--643.
\bibitem {Res87} Resnick, S., {\it Extreme Value Distributions, Regular
Variations, and Point Processes}, Springer, Berlin, 1987.
\bibitem
{Tad80} Tadikamalla, P.R. A look at the Burr and related
distributions. {\it International Statistical Review}, 48(1980)
337-344.
\bibitem {Urb56} Urbanik, K. On a problem concerning birth and
death processes. {\it Acta Math. Acad. Sci. Hungar.}, 7(1956),
99-106 (In Russian.)
\bibitem
{Wil94} Wilms, R. Fractional parts of random variables. Limit
theorems and infinite divisibility. Ph.D. Thesis, Technical
University of Eindhoven, Eindhooven, Holland, 1994.
\bibitem
{YanTso00} Yanev, G.P. and Tsokos, C.P. Family size order
statistics in branching processes with immigration. {\it Stoch.
Anal. Appl.} 18(2000), 4:655-670.
\bibitem {Zol54} Zolotarev, V.M. On a problem in the theory of
branching rpocesses. {\it Uspehi Matemat. Nauk} 9(1954), 147-156
(In Russian.)
\end{thebibliography}
\end{document}